\documentclass[10pt]{amsart}

\usepackage{amsmath,amssymb,amscd,enumerate,bbm}

\newcommand{\tmem}[1]{{\em #1\/}}
\newcommand{\tmop}[1]{\ensuremath{\operatorname{#1}}}
\newcommand{\tmstrong}[1]{\textbf{#1}}
\newcommand{\um}{-}
\newcommand{\tmname}[1]{\textsc{#1}}

\newcommand{\titl}{CLASSIFICATION PROBLEMS \\
AND MIRROR DUALITY}
\title{{\titl}}
\author{V. Golyshev}
\def\sec#1{\refstepcounter{section}\par\vspace{0.9cm}
\par\noindent
{\large\bf\arabic{section}. \ #1    } \addtocontents{toc} {\hbox
to\textwidth{\arabic{section}. #1 \dotfill \arabic{page}}}
\nopagebreak\par\vspace{0.3cm}}

\newcommand{\cal}{\mathcal}

\def\R{{\Bbb{R}}}
\def\C{{\Bbb{C}}}
\def\Q{{\Bbb{Q}}}
\def\Z{{\Bbb{Z}}}

\def\O{{\cal{O}}}
\def\E{{\cal{E}}}

\def\H{{\cal{H}}}
\def\D{{\cal{D}}}

\def\L{{\cal{L}}}

\def\Hom{{\rm H}{\rm o}{\rm m}}

\newcounter{pphcounter}[section]
\renewcommand{\thepphcounter}{\thesection.\arabic{pphcounter}}
\newcommand{\pph}{\bigskip \refstepcounter{pphcounter}
    \bf  \thepphcounter. \rm\quad}

\newcommand{\coro}{\bf Corollary. \rm}
\newcommand{\defi}{\bf Definition. \rm}
\newcommand{\rema}{\bf Remark. \rm}
\newcommand{\propo}{\bf Proposition. \rm}
\newcommand{\tensor}{\otimes}
\newcommand{\theo}{\bf Theorem. \rm}

\renewcommand{\proof}{\bigskip \bf Proof. \rm}

\renewcommand{\O}{{\mathcal{O}}}

\newcommand{\dimens}{{\mathrm{dim \;} }}

\def\A1{{{\Bbb{A}}^1}}
\def\P1{{{\Bbb{P}}^1}}

\newcommand{\Spec}{{\text{Spec }}}
\def\Anr1{{{\Bbb{A}} (n,r+1)}}

\def\Gm{{\bf G_m}}

\def\SL2{{\mathrm SL2}}

\renewcommand{\epsilon}{\varepsilon}

\def\conv{{*}}

\renewcommand{\P}{{\Bbb{P}}}
\renewcommand{\phi}{{\varphi}}

\newcommand{\eqdef}{\stackrel{def}{=}}

\newcommand{\lto}{\longrightarrow}

\newcommand{\mapto}{\mapsto}

\newcommand{\ind}{{\mbox{{$\mathrm{ind }$}}}}
\newcommand{\etaff}[4]{{{{$\frac{\mathbf{#1}^{#2}}{\mathbf{#3}^{#4}}$}}}}
\newcommand{\etaffff}[8]{{{{$\frac{\mathbf{#1}^{#2}\mathbf{#3}^{#4}}
{\mathbf{#5}^{#6}\mathbf{#7}^{#8}}$}}}}
\newcommand{\mystrut}{{\rule[-10pt]{0pt}{30pt}}}

\newcommand{\tns}{T_{\mathrm{{NS}}^\vee}}
\newcommand{\bbfC}{\C}

\newcommand{\inv}{{\mathrm{inv}}}
\newcommand{\FT}{{\mathrm{FT}}}
\newcommand{\symb}[2]{\left( \frac{#1}{#2}\right)}

\begin{document}

\begin{center}
\maketitle

\bigskip

\end{center}

\bigskip

In this paper we make precise, in the case of rank 1 Fano 3-folds,
the following program:
 {\em Given a classification problem in algebraic geometry,
use mirror duality to translate it into a problem in differential equations;
solve this problem and translate the result back into geometry. }

The paper is based on the notes of the lecture series the author gave
 at the University of Cambridge in 2003. It expands the announcement
\cite{Go2}, providing the background for and discussion of the modularity
conjecture.

We start with basic material on mirror symmetry for Fano varieties.
The quantum D-module
and the regularized quantum D-module are introduced in section 1. We state
the mirror symmetry conjecture for Fano varieties. We give
more conjectures implying, or implied by, the mirror symmetry
conjecture. We review the algebraic Mellin transform of Loeser and
Sabbah and define hypergeometric D-modules on tori.

In section 2 we consider Fano 3-folds of Picard rank 1 and review
Iskovskikh's classification into 17 algebraic deformation families.
We apply the basic setup to
Fano 3-folds to obtain the so called counting differential equations
of type D3. We introduce DN equations as generalization of these,
discuss their properties and take a brief look at their singularities.

In section 3, motivated by the Dolgachev-Nikulin-Pinkham picture of
mirror symmetry for K3 surfaces, we introduce
$(N,d)$-modular families; these are pencils of K3 surfaces whose
Picard-Fuchs equations are the counting D3 equations of rank 1 Fano
3-folds. The $(N,d)$-modular family is the pullback of the twisted
symmetric square of the universal elliptic curve over $X_0(N)^W$
to a cyclic covering of degree $d$.

Our mirror dual problem is stated in section 4: For which pairs
$(N,d)$ is it possible for the Picard-Fuchs equation of the corresponding
$(N,d)$-modular family to be of type D3? Through a detailed analysis
of singularities, we get a necessary condition on
$(N,d),$ bringing the list down to 17 possibilities.

Identifying certain odd Atkin-Lehner, weight 2, level
$N$ Eisenstein series (that appear in section 5)
with the sections of the bundle of relative differential 2-forms
in our modular family, we compute the corresponding Picard-Fuchs
equations and show them to be of type D3, recovering the matrix coefficients.

It turns out that
the pairs $(N,d)$ that we get are
exactly those for which there exists a rank 1 Fano 3-fold
of index $d$ and anticanonical degree $2d^2N.$
The Iskovskikh classification is revisited
in section 6. We sketch a proof that the matrices
we have recovered in section 5 via modular computations are, up to a
scalar shift, the counting matrices of the corresponding Fanos.

In section 7 we briefly discuss further classification problems to which
our approach may be applied.

\sec{Conjectures on Mirror Symmetry for Fano varieties.}

\pph \bf The big picture and the small picture. \rm There exist
two different approaches to differential systems built from
Gromov-Witten invariants of a variety.
The full Frobenius manifold underlies vector bundles
with connections
whose construction requires
knowledge of the \it big \rm quantum cohomology
and therefore of
the whole system of multiple-pointed correlators
(see chapter 2 in \cite{Ma}).
On the other hand, if
we are content to
restrict our study to the divisorial subdirection of the Frobenius manifold,
only the
\it small \rm quantum cohomology is needed. It is still a strong invariant of a
variety but it only requires knowledge
of the three-pointed correlators. For this reason, it is easier to compute.
The small quantum differential system has the additional advantage
of being representable as an algebraic D-module on the
``Neron-Severi-dual'' torus.

\pph \bf D-modules. \rm Given a smooth scheme
$X/\mathbb{C}$
we denote by
$D^{b, \text{holo}}$
the full subcategory of cohomologically bounded cohomologically holonomic
complexes of sheaves of left
$\mathcal{D}_X$-modules. For morphisms
$f\colon X \lto Y$ between smooth varieties, the ``six
operations'' exist and provide a convenient language for the
constructions that we are going to need. If
$f\colon X \lto Y$ is smooth of relative dimension
$d,$
then
$f_* K= Rf_*(K \tensor _{\mathcal{O}_X} \Omega^\bullet_{X/Y})[d]$, so that
$\H ^{i-d} (f_* K)= H^i_{DR} (X/Y,K)$ with its Gauss-Manin connection.

We will need the following notion of pullback: if
$M$ is a flat
$\mathcal{D}_Y$-module, then
$f^! M = f^+ M [ \dimens X - \dimens Y],$
where
$f^+M$
is the naive pullback of
$M$
as module with integrable connection.
 If
 $G$ is a separated smooth group scheme over
 $\C$
with group law
$\mu \colon G \times G \lto G,$
then the convolution of objects of
$D^{b, \text{holo}}(G)$ is defined by
$(K,L) \to K \conv L = \mu_* (K \boxtimes L)$
where
$K \boxtimes L$ is the external tensor product.

{\pph{ \bf Three-point correlators.}} Let $X$ be a Fano variety. Let $\tns$ be the
torus dual to the lattice $\tmop{NS}^{\vee}$(X). Define a trilinear functional
$\left\langle,, \right\rangle$ on the space $H ( X )$ setting

$$\left\langle \alpha, \beta, \gamma \right\rangle = \sum_{\chi \in
\tmop{NS}^{\vee} ( X )}
\left\langle \alpha, \beta, \gamma \right\rangle_\chi \cdot \chi $$
where
$\left\langle \alpha, \beta, \gamma \right\rangle _\chi$
is
``the expected number  of maps''
\label{dimension}
         \footnotemark
         \footnotetext{\parbox[t]{\textwidth}{
Technically, a Gromov-Witten invariant, \cite[VI-2.1]{Ma}. Let
$\overline{M}_n(X, \chi)$
denote the compactified moduli space of maps of rational curves of class
$\chi \in \mathrm{NS}^\vee$ with
$n$ marked points, and let
$\left[\overline{M}_n(X, \chi)\right]^\mathrm{virt}$
be its virtual fundamental class of virtual dimension
$\mathrm{vdim }\; \bar M_n(X, \chi)= \dimens X - \deg _{K_X} \chi +n -3.$
Let
$ev_i \colon \overline{M}_n(X, \chi) \lto X$ denote the evaluation map
at the
$i$-th marked point.
Then
$\langle \alpha, \beta, \gamma \rangle _\chi =
 ev_1^*(\alpha)\cdot  ev_2^*(\beta)\cdot  ev_3^*(\gamma)\cdot
 \left[\overline{M}_n(X, \chi)\right]^\mathrm{virt}$
if
$\mathrm{codim }\; \alpha +\mathrm{codim }\; \beta +\mathrm{codim }\; \gamma =
\mathrm{vdim }\; \overline{M}_3(X, \chi)$ and
$\langle \alpha, \beta, \gamma \rangle _\chi = 0$ otherwise.}}
         from $\mathbbm{P}^1$ to X in the cohomology class $\chi$ such
         that 0 maps into a general enough representative of
$\alpha,$ 1 maps into a representative of
$\beta$, $\infty$ maps into a representative of $\gamma$.
The functional $\left\langle,, \right\rangle$ takes values in
$\mathbbm{C}[ \chi]$.

{\pph{ \bf Algebra of quantum cohomology.}} {\tmstrong{{\tmstrong{}}}} Consider
the trivial vector bundle
$\mathcal{H} (X)$
over
$\tns$ with fiber $H(X).$  Extend the
Poincare pairing
$\left[\, , \right]$
to the vector space of its
sections
$H ( X ) \otimes \mathbbm{C} [ \chi]$. Raising an
index, we turn the trilinear form into a multiplication law on $H ( X ) \otimes
\mathbbm{C}[  \chi] :$

\begin{center}
  $\left[ \alpha \cdot \beta , \gamma \right] = \left\langle
  \alpha, \beta, \gamma \right\rangle$
\end{center}

{\pph{ \bf Connection (the $H^2 \um$direction).}} \label{connection} Identify
elements $f$ in the lattice $\tmop{NS} ( X )$ with invariant derivations
$\partial_f$ on $\tns$ by the rule
$${\partial_f ( \chi )} = f( \chi ) {\chi} $$
(In the left-hand side of the formula {\tmem{{\tmem{$\chi$}}}} is a function
on a torus; on the right, it is an element of the lattice $\tmop{NS} ( X )$
and as such is paired with $f$.)
Define a connection
$$
  \nabla_{\tns} \colon \Omega^0 ( \mathcal{H} ( X ) )
  \longrightarrow \Omega^1 ( \mathcal{H} ( X ) ) $$
in the vector bundle
$\mathcal{H} (X)$
by setting for any
constant section
$ \bar \alpha = \alpha \tensor 1 \in H ( X ) \otimes
\mathbbm{C}[  \chi]$
$$
 \left< \partial_f, \nabla_{\tns} \bar \alpha \right> = (f \tensor 1)\cdot \bar \alpha
$$
(the derivation $\partial_f$ is coupled with the vector-valued 1-form
$\nabla_{\tns} \bar \alpha$ in the LHS ).

\pph \theo The connection
$\nabla$
is flat.

\qed

This turns the space
$H ( X ) \otimes \mathbbm{C} [ \chi] $ into a
$\D =\D _{\tns}$-module. We will denote it by
$Q$ and call it the \emph{quantum D-module}.

{\pph{The mirror symmetry conjecture}} states that the solution
to the quantum D-module, convoluted with the canonical exponent, can
be represented as a period in some family of varieties, called a
(parametric) Landau-Ginzburg model. Let us make this more precise.

\pph \label{gma1}
\bf The exponent object on a one-dimensional torus. \rm Let
$\A1=\Spec \C[t],$
$\Gm=\Spec \C[t, t^{-1}],$ and let
$j\colon \Gm \lto \A1$ denote the corresponding open immersion.
Let $\partial=\frac{\partial}{\partial t}$
and
$D=D_t=t\frac{\partial}{\partial t}$
be the (invariant) derivations on
$\A1$ and
$\Gm,$ respectively.

The D-modules
$E = \mathcal{D}_\A1/\mathcal{D}_\A1 (\partial  -1)$,
and its restriction to
$\Gm$
$j^*E=\mathcal{D}_\Gm/\mathcal{D}_\Gm (D  -t),$
will be called the \emph{exponent objects.}

\bigskip

In general, the quantum
{\tmname{$\mathcal{D}$}}-module $Q$ is irregular. As such, it cannot possibly
be of geometric origin, that is, arise from a Gauss-Manin connection
of an algebraic family: Gauss-Manin connections are known
to be regular \cite{De}. In order to make a suitable geometricity
assertion one should pass to a regular object first.

{\pph{\bf Regularization.}}  Consider the inclusion
$\mathbbm{Z}K_X \rightarrow \tmop{NS}_X$. Dualizing twice, we have a morphism
of tori $\iota \colon \Gm \rightarrow T_{\tmop{NS}^{\vee}}$ (the canonical
torus map). Consider the exponent
object $j^* E$ on $\Gm$, and the pushforward $\iota_{\ast} (j^* E
)$. Define the \emph{regularized quantum object} as follows
\[ Q^{\tmop{reg}} = Q \ast \iota_{\ast} ( j^* E ) . \]

{\pph{ \bf The mirror symmetry conjecture.}} The object
$Q^{\tmop{reg}}$ is of geometric origin.

\pph \label{constituent}
This assertion in its strong interpretation means
that for any irreducible constituent
of
$Q^{\tmop{reg}}$
there exists a family
of varieties $\pi\colon \mathcal{E}\lto \tns$
such that
the restriction of that constituent
to some open subset
$U$
is isomorphic to a constituent of
$R^j \pi _* (\O).$
In practise (e.g. \ref{isolate--begin}--\ref{isolate--end})
we will forget about the trivial constituents
that may arise as a by-product of the convolution construction,
and deal only with the essential subquotient
of a single cohomology D-module of the regularized quantum object. We will call
it the \emph{regularized quantum D-module.}

Let
$\iota_x\colon \Gm \lto \iota (\Gm)x$ be an orbit of
$\Gm$ in
$\tns.$
The Gauss-Manin connection in the
Landau-Ginzburg
model with parameter
$x$
is then essentially the pullback
to
$\Gm$
of the regularized quantum D-module
with respect to
$\iota_x.$

\pph \label{strategy} Let us lay out broadly our classification strategy.
It is logical to start with the
Picard rank 1 case, as in this case
$\tns$ is one-dimensional and the regularized
quantum D-module is essentially a linear ordinary differential
equation with polynomial coefficients.

Assume we are interested
in finding all families of Fano varieties in a given class. (From our point
of view, a class comprises varieties with similar cohomology
structure. For instance, an interesting, if too narrow, class is that of
\emph{minimal Fanos} of a given dimension, i.e.
those whose non-trivial cohomology groups are just
$\Z$ in every even dimension. In the
class of \emph{almost minimal Fanos} we allow nontrivial
primitive cohomology in the middle dimension.) Assume that a variety
$X$ in the class is known, together with the values
$A_{X}=\{ a_{ij}(X)\}$ of the three-point correlators between two arbitrary
cohomology classes and the divisor class.
Compute the regularized quantum D-module and represent it as
$\mathcal{D}_\Gm/\mathcal{D}_\Gm L_{A_X}$ for
some
$L_{A_X} \in \mathcal{D}_\Gm$.
 We will say that
$L_{A_X}$  is  the counting differential operator
for $X$.
Doint the same construction starting with a matrix variable
$A=\{ a_{ij}\}$, we obtain a differential operator $L_A$ depending on
the set of parameters $\{a_{ij}\}$. (We will do this in detail for
almost minimal Fanos in \ref{dn}, getting what we call a differential
operator of type DN.) Thus, we can restate the
original classification problem as follows: determine which $L_A$ can
be counting differential operators $L_{A_X}$ of some Fano variety $X$.

\pph {\tmstrong{Identifying counting operators.}} What are the properties that
distinguish $L_{A_X}$'s as points in the affine space of all
$L_A$'s? As we have seen, the mirror symmetry conjecture asserts that the
$L_{A_X}$'s are of Picard-Fuchs type:
we expect that there exist  a pencil
$\pi\colon \E \lto \Gm$ defined over
$\Q$
and
$\omega$ a meromorphic section of a sheaf of relative differential
forms, such that a period
$\Phi$ of
$\omega$ satisfies
$L_{A_X} \Phi = 0.$
A believer in the mirror symmetry conjecture would therefore approach
the problem of identifying the possible $L_{A_X}$'s by first telling
which among all $L_A$'s are Picard-Fuchs. This will significantly narrow one's
search, as being Picard-Fuchs
is a very strong condition.

\pph \label{pf}
{\tmstrong{Identifying Picard-Fuchs operators}} among all $L_A$'s
apparently is not an algorithmic problem. The very first idea is to translate
(and that can be done algorithmically)
the   basic properties that an (irreducible) variation of Hodge structures
must have ---
regularity, polarizability, quasiunipotence of local monodromies ---  into
algebraic conditions on the coefficients of
the operator that represents it. One might hope
that these conditions cut out a variety of positive codimension
from the affine space of all $L_A$'s, thereby facilitating further search.
However, the hope is vain: a theorem proved recently by J. Stienstra and
myself asserts that a generic
DN equation is regular, polarizable and has quasiunipotent local monodromies
everywhere (see \ref{propertiesofdn}).

Algebraic  requirements  being met by virtue of the construction,
we have to shift the emphasis toward
non-algebrizable  conditions of analytic or of arithmetic nature
imposed by the PF property.

It is known that if a differential equation
$L_A \Phi =0$ with coefficients in
$\Q$  is of  Picard-Fuchs type,
then it is also

\begin{itemize}
  \item[({\bf H})] Hodge (that is, describes an abstract variation of $\Q$-Hodge
structures);
  \item[({\bf GN})] globally nilpotent
  \footnotemark
  \footnotetext{\parbox[t]{\textwidth}{
    We briefly recall what global nilpotence is. Let
    $\partial \xi = \xi M$ be an algebraic differential equation
    over
    $\mathbf{F}_p$. Consider
    $\partial  \partial \xi = \partial \xi M =\xi (M^2+M')$,
    $\partial  \partial \partial \xi$, etc. Then
    $\underbrace{\partial  \partial \dots \partial}_{p \text{ times}}
  \xi = C_p \xi$, for some matrix
    $C_p=C_p(M)$ which is called the
    $p$-curvature matrix. A differential equation
    $\partial \xi = \xi M$
    over
    $\Q$
    with
    $M$ having
    $p$-integral entries
    is said to be
    $p$-nilpotent
    if
    $C_p(M \mod p)$ is nilpotent. It is globally nilpotent if it is
    $p$-nilpotent for almost all
    $p$.}} in the sense of Dwork-Katz, see \cite{Dw}, \cite{Ka--NC}.
\end{itemize}

It is expected that, at least for small order
$r$ and degree
$d$, both ({\bf H}) and ({\bf GN}) are also sufficient conditions.
Unfortunately, there is no algorithmic way, given
$a_{ij},$
to verify that ({\bf H}) or ({\bf GN}) holds:
 in the former case, because of the fact that
({\bf H}) is, in particular, a condition on the global monodromy
which depends transcendentally on the coefficients of the
equation; in the latter case, because one does not know, given
$a_{ij},$ how to estimate the number of places
$(p)$ of $\Q$ where the nilpotence of the
$p$-curvature operator must be
verified in order to conclude that global nilpotence holds.

\pph
In order to state the hypergeometric pullback conjecture,
we will need some basic facts about hypergeometric
D-modules. Roughly, a D-module is hypergeometric
if the coefficients of the series expansion
of  its solution are products/quotients of the
gamma-function applied to values of nonhomogeneous linear forms in the degrees:
$\Phi = \sum u(n_1, \dots , n_p) t_1^{n_1}\; t_p^{n_p}$
with
$u(n_1, \dots , n_p)=
\prod c_i^{n_i} \prod_j \Gamma (l^{(j)}_i n_i - \sigma^{(j)})^{\gamma_j}.$
To put it precisely, one might use the language of algebraic Mellin transform,
introduced by Loeser and Sabbah (\cite{LS--EDF}).

%
Let
$\C [s]=\C [s_1,\dots, s_p]$ be the ring of polynomials in
$p$ variables and let $\C (s)$ be the corresponding fraction field.

\pph \label{relations} \defi A \emph{rational system of finite
  difference equations}
(FDE)  is a finite dimensional
$\C (s)$-vector space together with
$\C$-linear automorphisms
$\tau_1, \dots, \tau_p$ that commute with each other and satisfy the relations
$$ \tau_i s_j=s_j \tau _i  \text{  if  } i \ne j$$
$$ \tau_i s_i=(s_i + 1) \tau _i  \text{    } \forall i =1, \dots, p.$$

\pph If
$\frak M (s)$
and
$\frak M' (s)$
are rational systems of FDE, then so are
$\frak M (s)\tensor_ {\C (s)} \frak M' (s), \Hom_{\C (s)} (\frak M (s), \frak M' (s)).$
Therefore, the set of isomorphism classes of 1-dimensional
systems forms a group, which Sabbah and Loeser call the
\emph{hypergeometric group} and denote
$\mathcal{H}G(p).$

Denote by
$\L$ a subset of non-zero linear forms on
$\Q ^p$ with coprime integer coefficients
such that for all such forms
$L$
either
$L$
or
$- L$
is in
$\L. $
Let
$\Z ^ {[\L \times \C /\Z]}$ be
the subset of finitely supported functions
${\L \times \C /\Z} \lto \Z$
with the natural group structure.

\pph \propo \cite[1.1.4]{LS--EDF}. Let
$\sigma$ be a section of the projection
$\C \lto \C / \Z.$
Then, the map
$$(\C ^ *)^p \times \Z ^ {[\L \times \C /\Z]} \lto \mathcal{H}G(p)$$
that attaches to
$[(c_1, \dots, c_p), \gamma]$
the isomorphism class of the system satisfied by
$$(c_1)^{s_1} \dots (c_p)^{s_p} \prod_{L \in \L} \prod_{\alpha \in \C /
\Z} \Gamma(L(s)-\sigma(\alpha))^{\gamma_{L,\alpha}}$$
is a group isomorphism which
does not depend on the choice of
$\sigma.$

\pph Let
$T^p \simeq \Gm ^p$ be a complex torus of dimension
$p. $ Put
$D_i=t_i \frac{\partial}{\partial t_i}. $ Let
$\C [t,t^{-1}] \langle D \rangle$ denote the algebra of algebraic
differential operators on the torus (here
$t$ stands
for
$(t_1, \dots, t_p)$
and
$D$
for
$(D_1, \dots, D_p). $

The correspondence
$\tau_i=t_i$
and
$s_i=-D_i$
identifies this algebra with the algebra
$\C [s] \langle \tau, \tau ^{-1} \rangle$
of finite difference operators (which is the quotient
of the algebra freely generated by
$\C[s]$
and
$\C [\tau,\tau ^{-1}]$
by the relations
in Definition \ref{relations}).

If
$\mathcal{M}$ is a holonomic D-module on
$T^p,$ then its global sections form a
$\C [t,t^{-1}] \langle D \rangle$-module. The
\it algebraic Mellin transform \rm
$\frak M (\mathcal{M})$ of the D-module
$\mathcal{M}$ is this module of global sections
considered as
$\C [s]\langle \tau,\tau ^{-1}\rangle$-module. We say that
$\frak M (M)$
is a holonomic algebraic system of FDE if
$\mathcal{M}$ is holonomic.

\pph \bf The algebraic Mellin transform theorem. \cite[1.2.1]{LS--EDF}
\rm
Let
$\frak M $
be a holonomic algebraic system of FDE. Then
$\frak M (s) = \C (s) \tensor _{\C [s]} \frak M$
is a rational holonomic system of FDE.  Conversely, if
$\frak M (s)$
is a rational holonomic system of FDE, then for any
$\C [s]\langle \tau,\tau ^{-1} \rangle$-submodule
$\frak M \subset \frak M (s)$
such that
$\frak M (s) = \C (s) \tensor _{\C [s]} \frak M$
there exists
a holonomic algebraic system
$\frak M' \subset \frak M$
such that
$\frak M (s) = \C (s) \tensor _{\C [s]} \frak M'.$

\pph \propo \cite{LS--Ca} One has
 $\chi((\Gm)\sp p,\mathcal{M})=\dim\sb{\bbfC(s)}\frak M (\mathcal{M})(s)$.

\bigskip
Now we are ready to give a

\pph \defi
A D-module
$\mathcal{M}$
on
$T^p$
is said to be \emph{hypergeometric}
if
$\frak M (\mathcal{M})(s)$
has rank 1.

\bigskip

Every 1-dimensional $\bbfC(s)$- vector space
with invertible $\tau$-action contains a unique irreducible holonomic
$\C [s]\langle \tau,\tau ^{-1} \rangle$-module
and every such module of generic rank one is obtained in this way.

Passing back to the subject of quantum D-modules, we are finally set to state

\pph \bf The hypergeometric pullback conjecture. \rm
Let
$X$ be a Fano variety.
We conjecture that
for any constituent
$\mathcal{C}$
of the quantum D-module
$Q$
there exists a torus
$T_{\mathcal{C}},$ a morphism of
tori
$h_\mathcal{C}\colon \tns \lto T_\mathcal{C}$ and a hypergeometric D-module
$\mathcal{H_\mathcal{C}}
$
on
$T_\mathcal{C}$
such that
$\mathcal{C}$
is isomorphic to
a constituent of
the pullback
$h^!\mathcal{H_\mathcal{C}}$
on some open subset
$U$
of
$\tns.$

\pph \rema One can show
that
 the D-module
$Q$ is essentially  the restriction of the  ``extended first
structural connection'' onto the divisorial
direction the Frobenius manifold associated to
$X$   while
$Q^\mathrm{reg}$
corresponds to the ``second structural connection'', see
chapter 2 of \cite{Ma}.

\sec{The Iskovskikh classification and D3 equations.}

Let $X$ be a Fano 3-fold with one-dimensional Picard lattice, and
let $H=-K_X$ be the anticanonical divisor.
V.\ A.\ Iskovskikh classified all deformation families of these varieties
(see \cite{IP}). Recall that if $X$ is a smooth rank 1 Fano variety
and
$G \in H^2(X,\Z)$ is the positive generator
then the index of
$X$ is defined by
$H=(\ind X) G.$

\pph \label{classification}
\theo The possible pairs of invariants
$\displaystyle{(\frac{H^3}{2\, \mathrm{ind}^2 X}},\; \ind X)$ are
$$(1,1),(2,1),(3,1),(4,1),(5,1),(6,1),(7,1),(8,1),(9,1),(11,1),$$
$$(1,2),(2,2),(3,2),(4,2),(5,2),(3,3),(2,4).$$

\qed

To realize our strategy (\ref{strategy})
for rank one Fano 3-folds, one must first compute
the quantum D-module $Q.$

\pph \label{stable} \propo The subspace of algebraic classes
in the total cohomology
$H^\bullet (X)$ is stable under quantum multiplication by
$H.$ Therefore, the connection
$\nabla$ restricts to  the rank 4 subbundle of
$\mathcal{H}(X)$ generated by the algebraic classes.

\proof This follows easily from the ``dimension axiom'' (see the
formula in the footnote on
page \pageref{dimension}).

\qed

\pph \label{aij}
 We compute  this divisorial submodule explicitly, according to the definition.

Let us normalize
$a_{ij}$ so that
\bigskip
$$a_{ij}  = \begin{array}{l}
  \frac{1}{\text{deg} X} \cdot
 (j-i+1) \cdot \;
\text{the expected number of maps}\; \P ^1 \mapto X \\
\text{of degree}\; j-i+1\;
\text{that send 0 to the class of}\;
H^{3-i} \text{, and} \; \infty \; \text{to the class of} \; H^j.
\end{array}
$$

\bigskip

The degrees of the variety and of curves on it are considered
with respect
to
$H.$ Assume now for simplicity that
$X$ has index
$1$.

\bigskip

As always,
$\Gm = \Spec \C [t,t^{-1}]$ and
$D=t\frac{\partial}{\partial t}$. Let
$h^i$ be the constant sections of
$\mathcal{H}(X)$ that correspond to the classes
$H^i.$

\pph
\propo
The connection $\nabla$ is given by

$$
D (h^0, h^1, h^2, h^3)=(h^0, h^1, h^2, h^3)
\begin{pmatrix}
  a_{00}t & a_{01}t^2 & a_{02}t^3 & a_{03}t^4 \\
  1 & a_{11}t & a_{12}t^2& a_{13}t^3 \\
  0 & 1 &a_{22}t & a_{23}t^2 \\
  0 & 0 & 1 & a_{33}t
\end{pmatrix}
$$
\proof This follows from the definition.

\pph
\coro
Put
$$\widehat L_A= {\det}_{\text{right}} \left(D-\begin{pmatrix}
  a_{00}t & a_{01}t^2 & a_{02}t^3 & a_{03}t^4 \\
  1 & a_{11}t & a_{12}t^2& a_{13}t^3 \\
  0 & 1 &a_{22}t & a_{23}t^2 \\
  0 & 0 & 1 & a_{33}t
\end{pmatrix}\right).$$
where
${\det}_{\text{right}}$ means the ``right determinant'', i.e. the one that
expands as
\mbox{$\sum \text{\tiny element} \cdot \text{\tiny its algebraic
complement},$}
the
summation being over the
\it rightmost \rm column, and the algebraic complements being
themselves right determinants. Then
$h^0$ is annihilated by
$\widehat L$.

\proof This is a non-commutative version of Cayley-Hamilton.

\pph \coro
The quantum
D-module
$Q$
is isomorphic
to (a subquotient of)
$\mathcal{D}/\mathcal{D} \widehat{L}$.

\qed

Having thus computed $Q$, we proceed with regularization. We must convolute
$Q$ with the pushforward under the morphism
$\inv\colon x \mapsto 1/x$ of the exponent object $\D / (z \partial - z)\D$.
Convolution with the exponent of the inverse argument on a torus
is essentially
\footnotemark
the  Fourier (-Laplace) transform, as the following formula suggests:
$$(F(x) \conv (\frac{1}{x}e^{1/x}))(t)=\int F(y)\frac{y}{t}e^{y/t}\frac{dy}{y}=
 \frac{1}{t} (\FT(F))(\frac{1}{t}) $$
\footnotetext{Since everything is considered on
$\Gm$, the functions $\Phi$ and
$t \Phi$ are solutions to isomorphic D-modules. }

More precisely, one has:

\pph\bf The Fourier transform on
$\A1.$ \rm  \cite[2.10.0]{Ka--ESDE}  The \emph{Fourier transform}
of a differential operator
$L=\sum f_i(t) \partial^i \in \D_ \A1$
is defined by
$\FT(L)= \sum f_i(\partial)(-t)^i.$ The Fourier transform of the
left D-module
$M=\D_\A1/\D_\A1L$
is
$\FT (M) =\D_\A1/\D_\A1 \FT (L).$

\pph \label{isolate--begin} \propo \cite[5.2.3, 5.2.3.1]{Ka--ESDE}
Retain the notation of
\ref{gma1}.
 Then for any holonomic D-module
 $M$ on
 $\Gm$
we have
$$j^* \FT(j_* \, \inv _* (M)) \approx M \conv j^* E \quad
\text{and} \quad \inv_* j^* \FT(j_* M) \approx M \conv (\inv_* j^* E).$$

\pph \label{isolate--end}
The second formula shows that the convolution is in fact
a single D-module, though not in general an irreducible one.
We need to isolate
the essential subquotient, combing out the parasitic
 ones.
%

 Note that the operator
$\widehat L_A$ is divisible in
$\C [t, \partial]$ by
$t$ on the left (because the rightmost
column of the matrix is divisible
by $t$ on the left). Extend the D-module
$\D_\Gm/\D_\Gm t^{-1} \widehat L_A \approx  \D_\Gm/\D_\Gm  \widehat L_A$
naively to
$\A1$
as
$\D_\A1/\D_\A1 t^{-1} \widehat L_A.$
Do the Fourier transform. We get a D-module that corresponds to the
differential operator
$$\partial^{-1} {\det}_{\text{right}} \left(-D-1-\begin{pmatrix}
  a_{00}\partial & a_{01}\partial^2 & a_{02}\partial^3 & a_{03}\partial^4 \\
  1 & a_{11}\partial & a_{12}\partial^2& a_{13}\partial^3 \\
  0 & 1 &a_{22}\partial & a_{23}\partial^2 \\
  0 & 0 & 1 & a_{33}\partial
\end{pmatrix}\right).   \eqno{\mathrm{FT}}$$
Pass to the inverse: under
$\inv,$
$D$ is sent to
$-D$ and
$\partial$
to
$-t^2\partial.$ For further convenience we do two more things: shift
the differential operator by
$-1$ on the torus ($D$ goes to
$D$ and
$\partial$
to
$-\partial$) and multiply it by
$t^{}$ on the right. The result is then what we call a
\emph{counting differential
operator of type D3}. Abstracting our situation to any dimension and arbitrary
$\{a_{ij}\}$, we
introduce a

\pph  \label{dn} \defi
Let
$N$ be a positive integer. Let
$a_{ij} \in \Q, \;\; 0 \le i \le j \le N.$
Let
 $M$ be an $(N+1)  \times (N+1)$ matrix such that
for
$0 \le k,\, l \le N$:
$$
M_{kl}= \left\{
\begin{array}{ll}
0, & \text{ if } k > l+1,  \\
1, & \text{ if } k=l+1, \\
a_{kl} \cdot (Dt)^{l-k+1}, & \text{ if } k < l+1.
\end{array}
\right.
$$
We will also assume that the set
$a_{ij}$
is symmetric with respect to the SW-NE diagonal:
$a_{ij}=a_{N-j,N-i}$.

Put
$$\widetilde L = {\det}_{\text{right}} (D-M).$$
Since the rightmost column is divisible by
$D$ on the left,
the resulting operator
$\widetilde L$
is divisible by
$D$ on the left. Put
$$
\widetilde L = DL.
$$
The differential equation
$L\, \Phi(t)=0$ will be called a
 \it determinantal equation of order N, \rm or just a \emph{DN equation}.
Sometimes we write $\mathrm{DN}_{0,0}$ to signify that $0$ is a point
of maximally unipotent monodromy, and that the local expansion
$\Phi = c_0+c_1t + \dots$ of an analytic solution
$\Phi$ at $0$ starts with a nonzero constant term.
(One may have made other choices; for instance, the differential
operator marked $\FT$ above is of type $\mathrm{D3}_{\infty,1}$ in
this language.)

\pph \label{example} \bf Example. \rm A D3 equation expands as

\hspace{-12pt}\makebox{
$
\begin{array}{l}
 [{D}^{3}-t \left (2\,D+1\right )\left (a_{{00}}{D}^{2
}+a_{{11}}{D}^{2}+a_{{00}}D+a_{{11}}D+a_{{00}}\right )+ \\ + {t}^{2}
\left (D+1\right ) \left ({a_{{11}}}^{2}{D}^{2}+{a_{
{00}}}^{2}{D}^{2}+4\,a_{{11}}a_{{00}}{D}^{2}-a_{{12}}{D}^{2}-2\,a_
{{01}}{D}^{2}+ \right. \\
\phantom{aaaaaaaaaa}\left. + 8\,a_{{11}}a_{{00}}D-2\,a_{{12}}D+2\,{a_{{00}}}^{2}
D-4\,a_{{01}}D+2\,{a_{{11}}}^{2}D+6\,a_{{11}}a_{{00}}+{a_{{00}}}^
{2}-4\,a_{{01}}\right )- \\ - {t}^{3}\left (2\, D+3\right )\left (D+2\right
)\left (D+1\right )\left ({a_{{00}}}^{2}a_{{11}}+{
a_{{11}}}^{2}a_{{00}}-a_{{12}}a_{{00}}+a_{{02}}-a_{{11}}a_{{01}
}-a_{{01}}a_{{00}}\right )+ \\ + {t}^{4}\left (D+3 \right )\left
(D+2\right )\left (D+1\right )\left (-{a_{{00}}}^{2}a_{
{12}}+2\,a_{{02}}a_{{00}}+{a_{{00}}}^{2}{a_{{11}}}^{2}-a_{{03}}+
{a_{{01}}}^{2}-2\,a_{{01}}a_{{11}}a_{{00}}\right ) ] \, \Phi
(t) = 0
\end{array}
$
}

\pph
\defi
We say that two DN equations defined by sets
$a_{ij}$ and ${a_{ij}}^\prime$ \it are in the same class \rm if
there exists an
$a$ such that
$a_{ii}={a_{ii}}^\prime +a \; \text{for}\; i=0, \dots, N$
and
$a_{ij}={a_{ij}}^\prime \; \text{for}\; i \ne j$,
i.e.~if the matrices defined by
$a_{ij}$
and
${a_{ij}}^\prime$ differ by a scalar matrix.

Shifting the Fourier transformed
differential operator
$\FT$ on
$\A1$  corresponds exactly to shifting the DN matrix in its class.

\pph \defi We say that:

(i) a holonomic $\D$-module
$M$ is a variation of type
$DN$
if there exists a set of parameters
$A=\{a_{ij}\}$
such that
$\D / \D L_A \approx M.$
Here
$\approx$
denotes equivalence in the
category of
$\D$-modules up to
modules with  punctual support;

(ii) a constructible sheaf
$S$ is a variation of type
$DN$
if there exists a
$\D$-module
$M$
of type
$DN,$ such that
$H^{-1} (DR (M)) \approx S.$
Here
$\approx$
denotes equivalence in the
category of
constructible sheaves
up to sheaves with punctual support;
$DR$ is the Riemann-Hilbert correspondence functor.

\pph \label{propertiesofdn} \theo
(i)  A D-module  $\D / \D L$ of type DN is holonomic with
regular singularities;

(ii) it is self-adjoint;

(iii) the local monodromy around zero is
maximally unipotent (i.e. is conjugate to a Jordan block of size
$N$).

(iv) for a generic set
$A=\{a_{ij}\}$, the
D-module  $\D / \D L_A$ has
$N$ non-zero singularities. The local monodromies at those singularities
are symplectic (for
N even)
or orthogonal (for N odd) reflections, and the global monodromy is irreducible.

(v) the set  $A=\{a_{ij}\}$ can be recovered from the respective $L_A$:
if $A \ne A',$ then $L_A \ne L_{A'}.$
\bigskip

A proof can be found in a forthcoming paper by Jan Stienstra and myself.

\pph \defi We say that a DN variation
$M$ (resp.\ local system $S$)
is \emph{of geometric origin} if there exists a flat morphism
$\pi \colon \E \lto \Gm$
of relative dimension
$d$ such that $M$ (resp.\ $S$)
is isomorphic to a subquotient of the variation arising in its middle
relative cohomology ($R^0 \pi_*(\O)$, resp.\ $R^d \pi_*(\C)$)
up to a D-module (resp. a sheaf) with punctual support.

\pph \rema Recall that we had assumed in \ref{aij} that  the variety in
question had  index 1 before proceeding with the construction
of the counting differential operator. What happens in the higher index
cases? It turns out that
the definition \ref{dn} with the values of
$a_{ij}$ as defined in \ref{aij} is still valid, in the sense that it
yields a counting operator that corresponds to
the pullback
of the regularized quantum D-module
with respect to the anticanonical isogeny
$\Gm\stackrel{\mathrm{ind} X}{\lto}  \Gm$
(see \ref{constituent}.) We leave the proof to the reader. Use, for instance,
 the following property:

\pph \propo  \cite[5.1.9 1b]{Ka--ESDE} Let
$G$ be a smooth separated group scheme
of finite type,
$\phi \colon G \lto G$ a homomorphism. Then for any two objects
$K,L$ of
$D^{b, \text{holo}}(G)$
one has
$$ \phi^!((\phi_*K)\conv L) \approx K \conv (\phi^!L).$$

\pph \label{property}
In this language, the mirror symmetry conjecture for Fanos states:
the counting DN equations of almost minimal Fano $N$-folds
\footnote{One expects that the analog of Proposition \ref{stable}
holds in  even dimensions as well, so that all almost minimal $N$-folds
are controlled by DN's.}
are of geometric origin. In order to recover all counting DN equations
one should pose and then solve a mirror dual problem: find all
geometric DN equations that possess some special property. In general,
we do not know what that property is. However, in the D3 case we have
an additional insight: a counting D3 should come from an $(N, d)$-modular
family.

\pph \defi
A non-zero singularity of a D-module of type D3 is said to be:

(i) simple, if the local monodromy around that singularity is a
reflection (i.e. conjugate to the operator
$diag(-1,1,1)$);

(ii) complex, if it is not simple and is of determinant 1;

(iii) very complex, if it is not simple and is of determinant -1;

\sec{(N,\, d)-modular variation.}

\bf Warning. \rm   In this section
$N$ stands for level. This is not the
$N$ of the previous section, which denoted
the order of a differential operator.

{\pph{}} The quantum weak Lefschetz principle
implies that the fibers of the Landau-Ginzburg model of a Fano variety
are mirror dual to the sections of the anticanonical line bundle on
it. For rank 1 Fano 3-folds, these sections are rank 1 K3 surfaces.

The first picture of mirror symmetry for families of K3 surfaces arose as
an attempt to explain Arnold's strange duality. Let $L = 3U \oplus -2E_8$ be
the K3 lattice. For a wide class of primitive sublattices $M$ of
$L$ there is a unique decomposition
\[ M^{\bot} = U \oplus M^D, \]
so that there is a duality between
$M$ and
$M^D:$
\[ ( M^D )^{\bot} = U \oplus M. \] The Picard lattices of
mirror dual families of K3 surfaces are dual in this sense.
Therefore, it is natural
to expect that the dual Landau-Ginzburg model of a Fano 3-fold is
a family of K3 surfaces of Picard rank 19. We recall that a Kummer K3
is the minimal resolution of the quotient of an Abelian surfaces by
the canonical involution which sends $x$ to $-x$ in the group law.

\bigskip

The following construction was described in \cite{PS}, \cite{Go1}.

\pph \label{identification} Consider the modular curve
$X_0(N)$, and the
``universal elliptic curve'' over it.
Strictly speaking, the universal elliptic curve is a fibration not over
$X_0(N)$ but over a Galois cover with group $\Gamma$,
e.g.~$X(3N) - \{\text{cusps}\}$,
such that one can choose a $\Gamma$-form of the universal elliptic
curve; call it ``the'' universal elliptic curve and denote it by
$E_t$. Denote by $W$ the Atkin-Lehner involution of $X_0(N)$.
Consider the fibered product of
$E_t$ with the $N$-isogenous universal elliptic curve $E_t^W$ over
$X_0(N).$ We quotient out this relative abelian
surface $V_t$ by the canonical involution
$x \mapsto -x$ and then resolove to get a family of Kummer K3
surfaces. Let
$X_0(N)^\circ$ stand for
$X_0(N)-\{ \text{cusps} \} -\{ \text{elliptic points} \}$.
Denote by $H({V_{{\overline t_0}}})$ the cohomology of the generic
fiber of $V_t,$ that is, of the pullback of the family
$V_t$ to the universal cover of the base.

The monodromy representation
$$\psi\colon \pi _1(X_0(N)^\circ) \lto H^2(V_{{\overline t_0}})$$ is well
defined. We are going to compute
$\psi$ in
terms of the tautological projective representation
$$\phi \colon
\pi _1(X_0(N)^\circ) \lto PGL(H^1(E_{{\overline t_0}}))=PSL_2(\Z).$$

The monodromy that acts on
$H^1$ of the fiber of the universal elliptic curve is given by a lift of
$\phi$ to a linear representation
$$\bar \phi \colon \gamma  \mapto \begin{pmatrix}
  a & b \\
  c & d
\end{pmatrix}, \; c=0 \mod N.$$
Then, the monodromy that acts on
$H^1$ of the fiber of the isogenous curve is:
$$\bar \phi _N \colon \gamma \mapto \begin{pmatrix}
  d & -\frac{c}{N} \\
  -bN & a
\end{pmatrix}=
\begin{pmatrix}
0 & -\frac{1}{N} \\
1  & 0
\end{pmatrix}
\begin{pmatrix}
  a & b \\
  c & d
\end{pmatrix}
\begin{pmatrix}
0 & 1 \\
-N & 0
\end{pmatrix}
$$
where we have chosen symplectic bases
$\left< e_1,e_2 \right> ,\left< f_1,f_2 \right> $ of
$H^1(E_{{\overline t_0}}),H^1(E^W_{{\overline t_0}})$
such that the matrix of the isogeny $W$ in these bases
is
$$\begin{pmatrix}
0 & 1 \\
-N & 0
\end{pmatrix}.$$

The cohomology ring of the generic fiber
$V_{{\overline t_0}}$ of
our relative abelian surface is
$H(E_{{\overline t_0}}) \tensor H(E^W_{{\overline t_0}}).$
The vector subspace of algebraic classes in
$H^2(V_{{\overline t_0}})$ is generated by the pullbacks from the
factors and the graph of the isogeny:
$$
e_1 \wedge e_2 \tensor 1, 1 \tensor f_1 \wedge f_2,
- e_1 \tensor f_1 - N e_2 \tensor f_2.
$$
These classes are invariant under monodromy. The orthogonal
lattice of transcendental classes is generated by
$$ e_2 \tensor  f_1,
 e_1 \tensor f_1 -
 N e_2 \tensor f_2,
e_1 \tensor  f_2.$$

Identifying the $e$'s and $f$'s
with their pullbacks to the product, we write, abusing notation:
$$
e_2 \wedge f_1,  e_1 \wedge f_1-
 N e_2 \wedge f_2 ,  e_1 \wedge f_2.
$$
In this basis the monodromy representation is
$$
\psi\colon \gamma \mapto Sym_N^2 \, \phi (\gamma)=\begin{pmatrix}
  d^2 & 2cd & -c^2/N \\
  bd & bc+ad & -ac/N \\
  -Nb^2 & -2Nab &a^2
\end{pmatrix}
$$
(cf. \cite{PS}, \cite{Do}).

Let $\bar{\omega}$ be a meromorphic section of the sheaf of
relative holomorphic differential forms on the universal elliptic
curve. Identify $e_1, e_2 \text{ (resp. } f_1, f_2)$
with cohomology classes in the pullback of the universal elliptic curve
to the universal cover of the base. Denote by $\omega$ the pullback
of $\bar \omega$. Introduce a coordinate
$\tau$
on the universal cover by writing:
$$ [\omega] = \tau F e_1+Fe_2$$
(where $F$ is a function on the universal cover)
identifying it with the upper halfplane. The class $\omega ^W$ is then:
$$[\omega^W]=Ff_1-N\tau Ff_2.$$
Let $\omega$ and $\omega^W$ also denote, abusing notation,
the pullbacks of the respective forms to
$V_{\bar t_0}.$ Clearly,
$$[ \omega \wedge  \omega ^W] =
 F^2 e_2 \wedge f_1 +
 \tau F^2 (e_1  \wedge f_1
- N e_2 \wedge f_2) - \tau ^2 N F^2 e_1 \wedge f_2.
$$
Now, as
$ \omega$ is
$\Gamma_0 (N)$-equivariant,
$$\psi (\gamma)
\begin{pmatrix}
  F^2(\tau) \\ \tau F^2(\tau) \\-N \tau^2(F^2(\tau)
\end{pmatrix}
=
\begin{pmatrix}
  F^2(\gamma(\tau)) \\ \gamma(\tau)) F^2(\gamma(\tau)) \\-N \gamma(\tau))^2(F^2(\gamma(\tau))
\end{pmatrix},
$$
where
$\gamma(\tau)=\frac{a\tau +b}{c\tau +d}.$
This is equivalent to the identity
$$ F^2(\frac{a\tau +b}{c\tau +d})=(c \tau + d)^2 F^2(\tau).$$
Therefore, the period
$F^2$ in our family of abelian surfaces, as a function of
$\tau$, is a
$\Gamma_0(N)$-automorphic function of weight 2 on the upper halfplane.
Now, for any $\Gamma_0(N)$-automorphic function of weight 2 $G$,
the quotient
$\frac{G}{F^2}$
is
$\Gamma_0(N)$-invariant on the upper halfplane, hence a rational
function on
$X_0(N).$ This identifies
$G$
with a (meromorphic) section of the sheaf of
relative holomorphic 2-forms in our family.

Finally, delete the
$W$-invariant points from
$X_0(N)^\circ$
and let
${X_0(N)^W}^\circ$ be the quotient of the resulting curve by $W$.
The involution
$W$
extends to the fibration
$V_t$ in an obvious way, and yields a family
$V_t^W$ over
${X_0(N)^W}^\circ.$
The fundamental group
${X_0(N)^W}^\circ$
is generated by
$ \pi _1(X_0(N)^\circ)$
and
a loop
$\iota$ around
the point that is the image of a point
$s$ on the upper halfplane
stabilized by
$\begin{pmatrix}
0 & 1 \\
-N & 0
\end{pmatrix}.$
Extend
$\psi$ to
$\iota$, by
setting
$\psi (\iota)=\begin{pmatrix}
  0 & 0 & 1 \\
  0 & 1 & 0 \\
  1 & 0 & 0
\end{pmatrix}.$
The resulting representation is the monodromy representation of
the family
$V_t^W$ over
${X_0(N)^W}^\circ.$

If a relative holomorphic form
in the family
$V_t$
is a pullback from
$V_t^W,$
then, denoting its first period by $G$,
one has

$$\psi (\iota) \begin{pmatrix}
  G(\tau) \\ \tau G(\tau)  \\ -N \tau^2 G(\tau)
\end{pmatrix}
=
\begin{pmatrix}
  G(\frac{-1}{N \tau}) \\ \tau G(\frac{-1}{N \tau})  \\ -N \tau^2 G(\frac{-1}{N \tau})
\end{pmatrix},
$$
and
$G$ is odd Atkin-Lehner, as, by definition,
$$G^W(\tau)=G(\frac{-1}{N \tau})N^{-1}\tau ^{-2}.$$

Let now $N$ be a level such that the curve
$X_0(N)^W$ is rational.
We choose a coordinate
$T$ on it such that
$T=0$ at the image of the cusp
$(i \infty);$
(the inverse of a Conway-Norton uniformizer, see \ref{cn} below)
this defines an immersion of the torus
$\iota \colon X_0(N)^W \hookleftarrow \Gm^\prime = \Spec \C[T,T^{-1}].$
Let
$\Spec \C[t,t^{-1}] =\Gm \lto \Gm^\prime$ be the Kummer covering of
degree
$d,$ given by the homomorphism
$T \mapto t^d.$

The pullback of the variation described above
(that is, the pullback of the family itself, or
the monodromy representation, or
the
$D$-module, depending on the context)
to
$\Gm$
will be called the \emph{$(N,\, d)$-modular variation.}
Let us emphasize:
$(N,\, d)$-modular variations are variations on \emph{tori} ,
even if we speak of them as of variations on $\P ^1$,
as in the proof of Theorem \ref{maintheo} below.

In the case $N=1$  the construction is modified,
since the ``Atkin-Lehner involution''
$\begin{pmatrix}
0 & 1 \\
-N & 0
\end{pmatrix}$
acts trivially on $X_0(N)$.
In this case we work with the fibered product of
the ``universal elliptic curve'' over
$X_0(1)$ with its quadratic twist with respect to the degree two branched
covering ramified at the two elliptic points. In this case
the relative 2-form can no longer be identified
with a weight 2 level 1 modular function because of the sign
multiplier. However, squaring the corresponding period, we get
a bona fide modular function of weight 4 and level 1.

\sec{$(N,d)$-modular D3 equations: the necessary condition.}

\pph \bf Problem. \rm Find all pairs
$N,d$ such that the
$(N,d)$-modular variation described in the previous section
is of type $D3.$

\pph \label{maintheo}
\theo (Necessary condition). If an
$(N,d)$-modular variation is of type
D3, then the pair
$(N,d)$ belongs to the set

$$ \mathcal{M}=
\{ (1,1),(2,1),(3,1),(4,1),(5,1),(6,1),(7,1),(8,1),(9,1),(11,1),$$
$$(1,2),(2,2),(3,2),(4,2),(5,2),
(3,3),
(2,4) \}.$$

\bigskip
\proof We begin by noticing that no case with $d > 5$ is possible
as there would have to be at least 6 singularities.

\bf Case $\mathbf{d=1.}$ \rm

Assume $N \ne 1.$ We make the following remarks:

(1) \it All ramification points of the quotient map
$$\sigma \colon X_0(N)\lto X_0(N)^W $$
map to singularities of the $(N,1)$-modular variation. \rm
The corresponding local monodromy is projectively
(dual to) the symmetric square of the element in
$\Gamma_0 (N) + N$ that stabilizes this ramification point.
This element is elliptic or cuspidal, therefore its symmetric square
cannot be a scalar.

(2) \it Every elliptic point or a cusp point $p$ on $X_0(N)$ maps
to a complex or very complex point
$\sigma (p)$ on $X_0(N)^W.$ \rm If
$\sigma (p)$ were an apparent
singularity or a simple singularity, then the local monodromy
around $p$ would vanish,
which is precluded by the reason given above in (1).

(3) \it  The point
$s$ on the upper halfplane is neither elliptic nor a cusp. It goes to
a simple point on $X_0(N)^W.$ \rm We \sl defined \rm  the monodromy
$\iota (\text{image of } s )$ in the previous section to be a
reflection.

(4) \it  If a D3 equation is $(N,1)$-modular, then its set of non-zero
singularities consists of either 4 simple points,
or of 1 complex and 2 simple points, or of 1 very complex and
1 simple point. \rm The non-zero singularities of a D3
equation are inverse to roots of a polynomial of degree 4, as can
be seen  from the expansion in Example 1.3.
It has one simple singularity, according to (3).
Any singularity
of multiplicity 1 is simple. A singularity of multiplicity 3
is very complex because  the determinant must be $-1$ and it
cannot be simple (otherwise the global monodromy would be
generated by two reflections and therefore would be reducible).

(5) \it The genus $g$ of $X_0(N)$ is related to the numbers of elliptic points
$\nu_2$ and $\nu_3$ of order 2 and 3 on
$X_0(N)$ by the formula
$$g = 1 + \frac{N}{12} \prod_{p \mid N} (1+ p^{-1}) -
\frac{\nu_2}{4}-\frac{\nu_3}{3}-\frac{\nu_\infty}{2}.$$
\rm
This is Proposition 1.40 from \cite{Sh}.

These remarks show that
$g \le 1,$ (otherwise the variation would have at least 6 singularities
according to (1)); that if
$g=1$, then all of the singularities are simple (this is from (1) and (4)) and
$\nu_2=0, \; \nu_3=0, \nu_\infty=2$ so
$N=11$; and that if $g=0,$
then $N < 12$ (otherwise there would be too many singularities,
which would contradict (2) and (4)). The last argument also shows
that
$N \ne 10,$ as in this case
$\nu_2=2$ and
$\nu_\infty =4.$

\bf Case $\mathbf{d=2.}$ \rm

Again, assume $N \ne 1.$

(1) \it The genus $g$ of $X_0(N)$ is zero. \rm If it were greater than zero,
there would  be at least four singularities besides the one at
$0.$ Therefore, the $(N,2)$-modular variation would have at least
$7$ singularities.

(2) \it There may be no more than $3$ cusps
on $X_0(N).$
\rm
Assume there are at least 4 cusps on
$X_0(N).$
Consider the ramification points of the Atkin-Lehner involution.
One of them being
$s$,
the other is either a cusp or not a cusp.
In the former case we get at least three
cusps on
$X_0(N)^W.$ Pulling them back we get
at least
$4$ singularities of a D3
variation that are not simple, a contradiction. In the latter case,
we get at least two cusps and at least two other singularities of
the $(N,1)$-modular variation. Pulling them back to the
$(N,2)$-modular variation we get either:
\begin{itemize}
\item at least four simple points and two non-simple points, or:
\item at least two simple points and three non-simple points,
\end{itemize}
and in neither case can the resulting variation be of type D3.

(3) \it There may be no more than
$1$ order 3 elliptic point on
$X_0(N).$
\rm
Proof: same as above.

(4) \it There may be no more than
$7$ order 2 elliptic points on
$X_0(N).$
\rm These would give at least
$5$ singularities on
$X_0(N)^W$ and therefore at least
$7$ singularities on the pullback.

(5) \it The level
$N$ is smaller than
$48.$ \rm This bad but easy estimate
follows from the genus formula and the above remarks.

Having made these remarks, one proceeds (for instance) by inspecting
the values
$g, \nu_2, \nu_3, \nu_\infty$
for all levels
$N<48.$  One uses the formulas of \cite[Proposition 1.43]{Sh}:

$$
\nu_2 =\left\{
\begin{array}{ll}
0, & \text{ if } N=4k,  \\
\frac{1}{2}\prod_{p \mid N} (1+\symb{-1}{p}), & \text{ if } N=4k+2, \\
\prod_{p \mid N} (1+\symb{-1}{p}), & \text{ if } 2 \nmid N;
\end{array}
\right.
$$

$$
\nu_3 =\left\{
\begin{array}{ll}
0, & \text{ if } N=9k,  \\
\prod_{p \mid N} (1+\symb{-3}{p}), & \text{ if } 9  \nmid N;
\end{array}
\right.
$$

$$
\nu_\infty =\sum_{d \mid N} \phi(\mathrm{gcd}(d, N/d)).
$$
where $\phi (n)$ is as usual
the number of positive integers
not exceeding $n$ and relatively prime to $n.$

One thus finds that the only levels that satisfy the reqirements
above are
$N=2,3,4,5.$

\bf Case $\mathbf{d=3.}$ \rm

Pulling back under the degree 3 map dramatically multiplies singularities;
the analysis, which goes along the same rails, is this time much easier and
leaves one with the only possibility of a curve of genus 0 that has just
2 cusps, 1 order 3 elliptic point and no order 2 points, which corresponds
to level 3.

\bf Case $\mathbf{d=4.}$ \rm
Yet easier. The curve must be of genus 0 and have
2 cusps, 1 order 2 elliptic point and no order 3 points. The level is 2.

\qed

\sec{$(N,d)$-modular D3 equations: a sufficient condition.}

\pph \label{sufficient} \bf Theorem. \rm
\footnotemark
\footnotetext{\parbox[t]{\textwidth}{
Rather, a `fact', as its proof requires computations
too tedious to be done by hand. Note that, however, in the cases of
 complete intersections in projective spaces, the respective
$(N,1)$-modular local systems are rigid and can be identified
with the global monodromies of given D3 equations
by comparison of local monodromies, see \cite{Go1}.
}
}
 For all pairs $(N,d)$ in $\mathcal{M}$
the corresponding $(N,d)$-modular variation is of type
D3.

\proof
Assume for simplicity that $d=1$. Reshape the assertion this way:
for any pair
$(N,1)$ in $\mathcal{M}$
there is a period
$\Phi$  of a section of the line
bundle
$\pi _* \Omega ^2_{V_t/X_0(N)^W}$
in our
$(N,1)$-modular variation that satisfies, as a multivalued function on
$X_0(N)^W$,
a D3 equation with respect to a coordinate
$t$ on $X_0(N)^W.$ Given an expansion of
$\Phi(t)$ as a series in
$t,$ it is easy to find the differential equation that it satisfies.

To be more specific, recall that we chose a coordinate
$T$ on $X_0(N)^W$ in
\ref{identification} such that
$T=0$ at the image of the cusp
$(i \infty).$ The local monodromy at
$T=0$ of the cycles against which
our fibrewise 2-form is integrated
is conjugate to a unipotent Jordan block of
size $3$. Therefore, the analytic period
$\Phi=\Phi_0$ is well defined as the integral
against the monodromy-invariant cycle. In the same way,
the logarithmic period
$\Phi_1$, being the integral against a cycle in the second step
of the monodromy filtration, is well defined up to
an integral multiple of the analytic period. This defines
$\tau$ locally as
$\frac{\Phi_1}{\Phi_0}$,
and
$q$ as
$\exp ( 2 \pi i \frac{\Phi_1}{\Phi_0})$.

Now $q$ being a local coordinate around $0$, one can expand both
$\Phi$ and $T$ as $q$-series. Note that the expansion
of $T^{-1}$ is a
$q$-series that is uniquely defined up to a constant term. The
$q$-expansions of coordinates on
$X_0(N)^W$ appeared in a paper by Conway and Norton \cite{CN} and are called
Conway-Norton uniformizers. The table \ref{cn} of the uniformizers
for the levels that we need is taken from \cite{CN}.

Recall also that we have identified periods
$\Phi$ with odd Atkin-Lehner weight 2 level $N$ modular functions
in \ref{identification}.
Therefore, to prove our theorem explicitly one may: (1) produce a
$q$-expansion of such a
modular function
$\Phi$; (2)  fix the constant term  in the uniformizer
$T^{-1}$; (3) express
$q$ in
$T$; (4) expand
$\Phi$
in
$T$; (5) recover the differential equation that
$\Phi$ satisfies with respect to
$T$. If it is a D3 equation, we are done.

The same essentially  goes for the cases
$d=2,3,4$, except that the coordinate
on the Kummer covering is
$t=T^{1/d}$
and the local parameter is
$Q=q^{1/d}.$ The tables
\ref{table} contain the $Q$-expansions
of $\Phi$, the recovered
D3 matrices and the eta-expansions  of the
$I$-function that we introduce below.
The uniformizers, as we said, are next in the table
\ref{cn}. For level
$N$ and index
$d$ one should set the constant term of the uniformizer to
$a_{11}$ in the
$(N,1)$ matrix in tables \ref{table} (e.g.~take $c=744$ for level
$1$ and index $2$).

\bigskip
\bigskip
\bigskip

\pph \label{cn} \bf
The Conway-Norton uniformizers.  \rm
The constant term is denoted indiscriminately by $c$ below. We put
\mbox{$\mathbf{i}=q^{i/24} \prod (1-q^{in})$} in this table.

\bigskip

\hspace{-2pt}\makebox[520pt][r]{%
\begin{tabular}{|c|c|c|c|c|} \hline \mystrut
N=1 & N=2 & N=3 & N=4 & N=5 \\
\hline \mystrut
$j+ c$ &
\etaff{1}{24}{2}{24} + 4096 \etaff{2}{24}{1}{24}+ c &
\etaff{1}{12}{3}{12} + 729 \etaff{3}{12}{1}{12} + c&
\etaff{1}{8}{4}{8} + 256 \etaff{4}{8}{1}{8} + c&
\etaff{1}{6}{5}{6} + 125 \etaff{5}{6}{1}{6} + c\\
\hline
\hline \mystrut
N=6 & N=7 & N=8 & N=9 & N=11 \\
\hline \mystrut
\etaffff{1}{5}{3}{1}{2}{1}{6}{5} + 72 \etaffff{2}{1}{6}{5}{1}{5}{3}{1}  + c&
\etaff{1}{4}{7}{4} +  49 \etaff{7}{4}{1}{4} + c&
\etaffff{1}{4}{4}{2}{2}{2}{8}{4} + 32  \etaffff{2}{2}{8}{4}{1}{4}{4}{2} + c&
\etaff{1}{3}{9}{3} + 27  \etaff{9}{3}{1}{3} + c&
\etaffff{1}{2}{11}{2}{2}{2}{22}{2} + 16
\etaffff{2}{2}{22}{2}{1}{2}{11}{2} +
16  \etaffff{2}{4}{22}{4}{1}{4}{11}{4}+ c \\
\hline
\end{tabular}%
}

\bigskip

\pph \bf
Solution $\Phi.$ Weighted `Eisenstein series'
$E_2.$  \rm
In most of the cases, the form
$\Phi$  will be expressed as a finite linear combination of
``elementary Eisenstein series''
$$E_{2,i}(Q) \eqdef -\frac{1}{24}\,i \, (1 -  24 \sum_{n=1}^{\infty} \sigma (n)
Q^{in}).$$
A sequence $e_1,e_2,e_3, \dots$ determines the Eisenstein series
$$ \sum e_j \, E_{2,j}(Q).$$
We use notation
$\Phi = e_1 \cdot \mathbf{[1]} + e_2 \cdot \mathbf{[2]} + \dots $
in the third column of the  tables \ref{table}.

%

\pph \bf Non-uniqueness. \rm We have proved Theorem \ref{sufficient}
by producing \it some \rm
modular function
$\Phi$ and  \it some \rm  Conway-Norton uniformizer
$T^{-1}$
of level
$N$
such that
$\Phi$ expanded in
$T$ satisfies a D3 equation.
Is the pair
$\Phi, \; T^{-1}$
that we have produced
determined
by this condition uniquely?  The answer is in general no, even if
$\Phi$ is known to be an Eisenstein series:
at certain composite levels the space spanned by Eisenstein series
has dimension higher than $1$, and it is possible to find two
different Eisenstein series and two uniformizers (that differ by a
constant term) such that the respective expansions give rise to
different D3 matrices.


The extra piece
that we use to characterize the matrices
and the solutions
$\Phi$
in the tables
\ref{table}  uniquely
is:

\pph \bf The miraculous eta-product formula. \rm
Define
$I=\Phi \cdot {t}^{d \, \frac{N+1}{12}}.$
Let \mbox{$ H_j (Q)=Q^{j/24} \prod (1-Q^{jn}).$}
It turns out that
$I$ expands
as a finite product of series of the form
$\prod H_j^{h_j}(Q)$ in a remarkably uniform way:
$$  I=H_d(Q)^2 H_{Nd}(Q)^2. $$
We reflect this phenomenon in the
fourth
column of the  tables \ref{table}.
The notation used is
$I = \mathbf{1}^{h_1} \mathbf{2}^{h_2} \cdot \dots \; .$
No intrinsic explanation of the eta-product formula is known to the author.

\pph \label{table} {\bf Level, matrix, solution,
$I$-function.}

\mystrut $d=1$

\vbox{
\hspace{0pt}\makebox[520pt][l]{
\begin{tabular}{|c|c|c|c|}
\hline \mystrut
 N
 & $a_{ij}$
 & $\Phi \; \; (\mathbf{[j]}=-\frac{1}{24}\,E_{2,j}(Q)) $
 & $I  \; \; (\mathbf{j}=H_j(Q))$\\ \hline
 1
 &$ \begin{array}{cccc}120 & 137520 & 119681280 & 21690374400 \\ 0 & 744 & 650016 & 119681280 \\ 0 & 0 & 744 & 137520 \\ 0 & 0 & 0 & 120 \end{array} $
&$\sqrt{E_4(q)}$
& $ \mathbf{1} ^2 \mathbf{1} ^2 $\\  \hline
 2
 &$ \begin{array}{cccc}24 & 3888 & 504576 & 18323712 \\ 0 & 104 & 13600 & 504576 \\ 0 & 0 & 104 & 3888 \\ 0 & 0 & 0 & 24 \end{array} $
& $ +24\cdot \mathbf{[1]}-24\cdot \mathbf{[2]} $
& $ \mathbf{1} ^2 \mathbf{2} ^2 $\\  \hline
 3
 &$ \begin{array}{cccc}12 & 792 & 43632 & 793152 \\ 0 & 42 & 2340 & 43632 \\ 0 & 0 & 42 & 792 \\ 0 & 0 & 0 & 12 \end{array} $
& $ +12\cdot \mathbf{[1]}-12\cdot \mathbf{[3]} $
& $ \mathbf{1} ^2 \mathbf{3} ^2 $\\  \hline
 4
 &$ \begin{array}{cccc}8 & 304 & 9984 & 121088 \\ 0 & 24 & 800 & 9984 \\ 0 & 0 & 24 & 304 \\ 0 & 0 & 0 & 8 \end{array} $
& $ +8\cdot \mathbf{[1]}-8\cdot \mathbf{[4]} $
& $ \mathbf{1} ^2 \mathbf{4} ^2 $\\  \hline
 5
 &$ \begin{array}{cccc}6 & 156 & 3600 & 33120 \\ 0 & 16 & 380 & 3600 \\ 0 & 0 & 16 & 156 \\ 0 & 0 & 0 & 6 \end{array} $
& $ +6\cdot \mathbf{[1]}-6\cdot \mathbf{[5]} $
& $ \mathbf{1} ^2 \mathbf{5} ^2 $\\  \hline
 6
 &$ \begin{array}{cccc}5 & 96 & 1692 & 12816 \\ 0 & 12 & 216 & 1692 \\ 0 & 0 & 12 & 96 \\ 0 & 0 & 0 & 5 \end{array} $
& $ +5\cdot \mathbf{[1]}-1\cdot \mathbf{[2]}+1\cdot \mathbf{[3]}-5\cdot \mathbf{[6]} $
& $ \mathbf{1} ^2 \mathbf{6} ^2 $\\  \hline
 7
 &$ \begin{array}{cccc}4 & 64 & 924 & 5936 \\ 0 & 9 & 140 & 924 \\ 0 & 0 & 9 & 64 \\ 0 & 0 & 0 & 4 \end{array} $
& $ +4\cdot \mathbf{[1]}-4\cdot \mathbf{[7]} $
& $ \mathbf{1} ^2 \mathbf{7} ^2 $\\  \hline
8
 &$ \begin{array}{cccc}4 & 48 & 576 & 3328 \\ 0 & 8 & 96 & 576 \\ 0 & 0 & 8 & 48 \\ 0 & 0 & 0 & 4 \end{array} $
& $ +4\cdot \mathbf{[1]}-2\cdot \mathbf{[2]}+2\cdot \mathbf{[4]}-4\cdot \mathbf{[8]} $
& $ \mathbf{1} ^2 \mathbf{8} ^2 $\\  \hline
 9
 &$ \begin{array}{cccc}3 & 36 & 378 & 1944 \\ 0 & 6 & 72 & 378 \\ 0 & 0 & 6 & 36 \\ 0 & 0 & 0 & 3 \end{array} $
& $ +3\cdot \mathbf{[1]}-3\cdot \mathbf{[9]} $
& $ \mathbf{1} ^2 \mathbf{9} ^2 $\\  \hline
 11
 &$ \begin{array}{cccc}12/5 & 24 & 198 & 880 \\ 0 & 22/5 & 44 & 198 \\ 0 & 0 & 22/5 & 24 \\ 0 & 0 & 0 & 12/5 \end{array} $
& $ +12/5 \cdot \mathbf{[1]}-12/5 \cdot \mathbf{[11]} $
& $ \mathbf{1} ^2 \mathbf{11} ^2 $\\  \hline
\end{tabular}
}
}
%

\newpage
\mystrut $d=2$

\hspace{-12pt}\begin{tabular}{|c|c|c|c|}
\hline \mystrut
N
& $a_{ij}$
& $\Phi
$
& $I$ \\ \hline
 1
 &$ \begin{array}{cccc}0 & 240 & 0 & 57600 \\ 0 & 0 & 1248 & 0 \\ 0 & 0 & 0 & 240 \\ 0 & 0 & 0 & 0 \end{array} $
&$\sqrt{E_4(Q^2)\mathstrut}$
& $ \mathbf{2} ^2 \mathbf{2} ^2 $\\  \hline
 2
 &$ \begin{array}{cccc}0 & 48 & 0 & 2304 \\ 0 & 0 & 160 & 0 \\ 0 & 0 & 0 & 48 \\ 0 & 0 & 0 & 0 \end{array} $
& $ +12\cdot \mathbf{[2]}-12\cdot \mathbf{[4]} $
& $ \mathbf{2} ^2 \mathbf{4} ^2 $\\  \hline
 3
 &$ \begin{array}{cccc}0 & 24 & 0 & 576 \\ 0 & 0 & 60 & 0 \\ 0 & 0 & 0 & 24 \\ 0 & 0 & 0 & 0 \end{array} $
& $ +6\cdot \mathbf{[2]}-6\cdot \mathbf{[6]} $
& $ \mathbf{2} ^2 \mathbf{6} ^2 $\\  \hline
 4
 &$ \begin{array}{cccc}0 & 16 & 0 & 256 \\ 0 & 0 & 32 & 0 \\ 0 & 0 & 0 & 16 \\ 0 & 0 & 0 & 0 \end{array} $
& $ +4\cdot \mathbf{[2]}-4\cdot \mathbf{[8]} $
& $ \mathbf{2} ^2 \mathbf{8} ^2 $\\  \hline
 5
 &$ \begin{array}{cccc}0 & 12 & 0 & 160 \\ 0 & 0 & 20 & 0 \\ 0 & 0 & 0 & 12 \\ 0 & 0 & 0 & 0 \end{array} $
& $ +3\cdot \mathbf{[2]}-3\cdot \mathbf{[10]} $
& $ \mathbf{2} ^2 \mathbf{10} ^2 $\\  \hline
\end{tabular}

\mystrut $d=3$

\begin{tabular}{|c|c|c|c|}
\hline \mystrut
N
& $a_{ij}
$ & $\Phi$
& $I$ \\ \hline
 3
 &$ \begin{array}{cccc}0 & 0 & 54 & 0 \\ 0 & 0 & 0 & 54 \\ 0 & 0 & 0 & 0 \\ 0 & 0 & 0 & 0 \end{array} $
& $ +4\cdot \mathbf{[3]}-4\cdot \mathbf{[9]} $
& $ \mathbf{3} ^2 \mathbf{9} ^2 $\\  \hline
\end{tabular}

\bigskip

\mystrut $d=4$

\begin{tabular}{|c|c|c|c|}
\hline \mystrut
N
& $a_{ij}$
& $\Phi$
& $I$ \\ \hline
 2
 &$ \begin{array}{cccc}0 & 0 & 0 & 256 \\ 0 & 0 & 0 & 0 \\ 0 & 0 & 0 & 0 \\ 0 & 0 & 0 & 0 \end{array} $
& $ +6\cdot \mathbf{[4]}-6\cdot \mathbf{[8]} $
& $ \mathbf{4} ^2 \mathbf{8} ^2 $\\  \hline
\end{tabular}

\pph \label{equal} \rema The functions
$\Phi$
and
$I$ that describe in these tables the cases
with the same level
$N,$ but with different
$d,$ are equal.

\newpage

\pph {\bf Respective differential operators.}

\mystrut $d=1$

\begin{tabular}{l|l}
   1
& ${D}^{3}-24\,t\left (1+2\,D\right )
\left (6\,D+5\right )\left (6\,D+1
\right )
$
                                  \\\hline     2
& ${D}^{3}-8\,t\left (1+2\,D\right )
\left (4\,D+3\right )\left (4\,D+1
\right )
$                                    \\\hline
                                    3
& ${D}^{3}-6\,t\left (1+2\,D\right )
\left (3\,D+2\right )\left (3\,D+1
\right )
 $

                                      \\\hline 4
&
${D}^{3}-8\,t\left (1+2\,D\right )^{3}

$

                                     \\\hline  5
& ${D}^{3}-2\,t\left (1+2\,D\right )
\left (11\,{D}^{2}+11\,D+3\right )-4\,
{t}^{2}\left (D+1\right )\left (2\,D+3
\right )\left (1+2\,D\right )

$

                                      \\\hline 6
&${D}^{3}-t\left (1+2\,D\right )\left (
17\,{D}^{2}+17\,D+5\right )+{t}^{2}
\left (D+1\right )^{3}

$

                                    \\\hline   7
&${D}^{3}-t\left (1+2\,D\right )\left (
13\,{D}^{2}+13\,D+4\right )-3\,{t}^{2}
\left (D+1\right )\left (3\,D+4\right
)\left (3\,D+2\right )

$

                                     \\\hline  8
&${D}^{3}-4\,t\left (1+2\,D\right )
\left (3\,{D}^{2}+3\,D+1\right )+16\,{
t}^{2}\left (D+1\right )^{3}

$

                                      \\\hline 9
&${D}^{3}-3\,t\left (1+2\,D\right )
\left (3\,{D}^{2}+3\,D+1\right )-27\,{
t}^{2}\left (D+1\right )^{3}

$

 \\\hline 11
&$ {D}^{3}-2/5\,t\left (2\,D+1\right )\left (17\,{D}^{2
}+17\,D+6\right )-{\frac {56}{25}}\,{t}^{2}\left (D+
1\right )\left (11\,{D}^{2}+22\,D+12\right )-$\\ & $\mystrut -{\frac {126}{125}}\,{t}^{
3}\left (2\,D+3\right )\left (D+2\right )\left (D+1
\right )-{\frac {1504}{625}}\,{t}^{4}\left (D+3
\right )\left (D+2\right )\left (D+1\right )

$ \\\hline

\end{tabular}

\mystrut $d=2$

\begin{tabular}{l|l}
   1
& ${D}^{3}-192\,{t}^{2}\left (3\,D+5
\right )\left (3\,D
+1\right )\left (D
+1\right )
$
                                  \\\hline     2
& ${D
}^{3}-64\,{t}^{2}\left (2\,D
+3
\right )\left (2\,D
+1\right )\left (D
+1\right )
$                                    \\\hline
                                    3
& ${D
}^{3}-12\,{t}^{2}\left (3\,D
+2
\right )\left (3\,D
+4\right )\left (D+1\right )
 $

                                      \\\hline 4
&
${D
}^{3}-64\,{t}^{2}
\left (D
+1\right )^{3}

$

                                     \\\hline  5
& ${D
}^{3}-4\,{t}^{2}
\left (D+1\right )\left (11\,{D}^{2}+
22\,D+12\right )-16\,{t}^{4}\left (D+3\right )\left (D+2
\right )\left (D+1\right )

$

                                      \\\hline

\end{tabular}

\bigskip

\bigskip

\mystrut $d=3$

\begin{tabular}{l|l}

                                    3
& $ {D}^{3}-54\,{t}^{3}\left (2\,D+3\right )\left (D+2
\right )\left (D+1\right )
$

                                      \\\hline

\end{tabular}

\bigskip

\bigskip

\mystrut $d=4$

\begin{tabular}{l|l}
    2
& ${D}^{3}-256\,{t}^{4}\left (D+3\right )\left (D+2
\right )\left (D+1\right )
$                                    \\\hline

\end{tabular}

\bigskip

\sec{A conjecture on counting matrices. The Iskovskikh classification
  revisited.}

\pph {\bf Corollary of theorems \ref{classification} and \ref{maintheo}.}
The $d$-Kummer pullback of the Picard-Fuchs equation
of the twisted symmetric square of the universal elliptic curve over
$X_0(N)^W$ is of type D3 if and only if there exists a family of rank
1 Fano 3-folds  of index $d$ and anticanonical degree  $2d^2N$.

\pph \label{conjecture}\bf Modularity conjecture. \rm The counting matrix of
a generic Fano 3-fold in the Iskovskikh family with parameters $(N,d)$
is in the same class as the corresponding matrix in tables \ref{table}.

\bigskip

More concretely, the conjecture states that the matrix $a_{ij}$ of normalized
Gromov-Witten invariants of a Fano 3-fold with invariants
$(N,d)$ can be obtained in the following uniform way. Let $T=T(q)$ be
the inverse of the suitable Conway-Norton uniformizer on
$X_0(N)$ (that is, the one with the ``right'' constant term). Consider
$$\Phi = (q^{1/24} \prod (1-q^{n})q^{N/24} \prod
(1-q^{Nn}))^2T^{-\frac{N+1}{12}}.$$
Then $\Phi$ satisfies a D3 equation with respect to
$t=T^{\frac{1}{d}}$. Recover the matrix of
$a_{ij}$ that corresponds to this equation (e.g. by looking at the expansion
in Example \ref{example}), and normalize it by subtracting $a_{00}I$.

This uniform description is somewhat unexpected, since it does not
have an obvious translation in terms of the geometry of Fano 3-folds.
Let us now take a more detailed view at the Iskovskikh classification,
according to the index and the degree.

\pph \bf The Iskovskikh classification revisited. \rm
\bigskip

\mystrut $d=1$

\begin{tabular}{l|l}

1  & hypersurface of degree $ 6$ in $\P (1,1,1,1,3) $ \\ \hline
2  & quartic  in $\P ^4 $ \\  \hline
3  & complete intersection of a quadric and a cubic in $ \P ^5$ \\  \hline
4  & complete intersection of 3 quadrics in $\P ^6$ \\  \hline
5  & a section of the Grassmannian $G(2,5)$ by a quadric and a codimension 2 plane \\  \hline
6  & a section of the orthogonal Grassmannian $O(5,10)$ by a codimension 7 plane \\  \hline
7  & a section of the Grassmannian $G(2,6)$ by a codimension 5 plane \\  \hline
8  & a section of the lagrangian Grassmannian $L(3,6)$ by a codimension 3 plane \\  \hline
9  & a section of $G_2/P$ by a codimension 2 plane \\ \hline
11 & variety $V_{22}$ \\ \hline

\end{tabular}

\bigskip

\bigskip

\mystrut $d=2$

\begin{tabular}{l|l}

1  & hypersurface of degree $ 6$ in $\P (1,1,1,2,3) $ \\ \hline
2  & hypersurface of degree $4 $ in $\P (1,1,1,1,2) $ \\ \hline
3  & a cubic in $ \P ^4$ \\  \hline
4  & complete intersection of 2 quadrics in $\P ^6$ \\  \hline
5  & a section of the Grassmannian $G(2,5)$ by a codimension 3 plane \\  \hline

\end{tabular}

\bigskip

\bigskip

\mystrut $d=3$

\begin{tabular}{l|l}

3  &  quadric  in $ \P ^4$ \\  \hline

\end{tabular}

\bigskip
\bigskip

\mystrut $d=4$

\begin{tabular}{l|l}

2  &   $ \P ^3$ \\  \hline

\end{tabular}

\bigskip

\pph \rema The description of families
$(6,1),(8,1),(9,1)$ as hyperplane sections in Grassmannians is due to Sh. Mukai
\cite{Mu}.

\pph \bf How can one prove the modularity conjecture?

\rm The uniformity of the assertion calls for a uniform proof, but I
do not know how such a proof might work.

The only way I know how to prove the conjecture is to explicitly calculate
the quantum cohomology of Fano 3-folds on a case by case basis.

Kuznetsov calculated the quantum cohomology of $V_{22}$.
All other cases are complete intersections in weighted projective
spaces or Grassmannians of simple Lie groups.

For complete intersections in usual projective space, Givental's
result allows to compute the D3 equations and the result agrees with
the conjecture. Przyjalkowski \cite{Pr} has recently extended
Givental's result to the cases of smooth complete intersections in
weighted projective spaces and established the predictions in the cases
$(N,d) \in \{(1,1),(1,2),(2,2)\}$.

In the remaining cases we use the quantum Lefschetz principle to
reduce the computation of the quantum D-module of a hyperplane section
to that of the ambient variety.

\pph \theo (Quantum Lefschetz hyperplane section theorem,
Coates-Givental-Lee-Gathmann, see e.g. \cite{Ga}.)
Let
$Y$ be a section of a very ample line bundle
$\L$ on $X$. We assume that both varieties are of Picard rank 1. Let
$\iota_{\widehat \L}\colon \Gm \rightarrow T_{\tmop{NS}^{\vee}}$
be the morphism of tori double dual to the
map
$\Z[\widehat \L] \lto NS_X.$ For
$\lambda \in \C^*$ let
$[\lambda] \colon \Gm \lto \Gm$ be the corresponding translation. For
$\alpha \in \C^*$ let
$[\alpha] \colon \A1 \lto \A1$ be the corresponding multiplication
morphism.
Then the quantum D-modules are related as follows:

(i) if the index of  $Y > 1$,
there exists
$\lambda$  in
$\C ^*$ such that
$$
Q_Y\;
\text{is a subquotient of}\;
[\lambda]^* (Q_X \conv {\iota_{\widehat \L}}_*(j^* E))
$$

(ii) if the index of  $Y = 1,$ there exist
$\lambda$ and
$\alpha$ in
$\C ^*$ such that
$$Q_Y\;
\text{is a subquotient of}\;
[\lambda]^* (Q_X \conv {\iota_{\widehat \L}}_*(j^* E)) \tensor j^* ([\alpha]_* E).
$$

\pph The quantum cohomology of ordinary, orthogonal and lagrangian
Grassmannians is known
(Givental-Kim-Siebert-Tian-Peterson-Kresch-Tamvakis).
Przyjalkowski calculated the quantum Lefschetz reduction for the cases
$(5,1),(7,1)$, confirming the conjecture.

Note that we do not need the whole cohomology structure: we just need
to know quantum multiplication by the divisor classes, and this can be
computed using Peterson's quantum Chevalley formula  \cite{FW}.
I calculated the quantum Lefschetz reduction for the cases
$(6,1),(8,1),(9,1)$ and the results again agreed with the ones
predicted by the conjecture.

To our knowledge, quantum multiplication by the divisor class on
$V_5$ (case $(5,2)$) was first computed by Beauville \cite{Bea}. We refer
the reader to \cite{BM} which makes use of Beauville's and Kuznetsov's results.

To summarize, we have checked the conjecture in all 17 cases by a case
by case analysis. This proof, however, does not explain why the
conjecture is true. A more uniform approach,
yet to be discovered, would presumably start from the embedded K3 rather than
the ambient space.

\pph  \rema If the conjecture is true, then there
is a mysterious relation between varieties of different index,
as implied by Remark \ref{equal}.

\sec{What next?}

\pph \bf Classification of smooth rank 1 Fano 4-folds.  \rm This is an
open question. For a variety of index
$\ge 2$ one can pass to the
hyperplane section (which has to be a Fano 3-fold) and thus
reduce the problem to lower dimension. On the other hand, the
classification of index one Fano 4-folds seems to be beyond reach
of today's geometric methods. Our program, if  carried  out in
this case, would suggest a blueprint of a future classification.

As a first step one must show that rank 1 Fano 4-folds do give rise
to equations of type D4. The dimension argument that we used in
\ref{stable} to show that the subspace
generated by
$H^n, \; n=0, \dots ,\dimens X$, is stable under quantum multiplication
by
$H$ no longer works. Still, the assertion is true in dimension 4.
The next step is to classify counting D4 equations.
Unlike D2 and D3 variations, whose differential Galois group is
$\mathrm{Sl}_2=\mathrm{Sp}_2=\mathrm{So}_3$,
a variation of type D4 is controlled by
$\mathrm{Sp}_4,$ and has in general no chance of being modular.
Thus, as we remarked in \ref{property}, in the D4 case we lack the
consequences of modularity that
enabled us first to state the correct mirror dual problem, and then effectively
to handle it in the D3 case.

With no idea of what the mirror dual problem might be, one can still
rely on the basic conjectures of chapter 1 to
compose a list of candidate D4 equations. If the list is not too long
and it contains all D4 equations, the problem is reduced
to weeding out the extra non-counting D4 equations that
have sneaked into the list.

Which D4 equations are of Picard-Fuchs type? Of the approaches
that we discuss in \ref{pf}, establishing the $\Q$-Hodge
or even the $\R$-Hodge property of a differential
equation, given its coefficients, seems hopeless. On the other hand,
a necessary, though not sufficient, condition for global nilpotence
is that the $p$-curvatures are nilpotent for sufficiently many prime
$p$. In principle   \footnotemark
\footnotetext{But not in practice. The generic D4
depends on 9 parameters, and the computation involved
needs unrealistic resources.}, one needs to guess
the upper bound $h_{\max}$ of the height ($h= \max (p,q)$
for $p/q \in \Q$ in lowest terms)
of possible Gromov-Witten invariants $a_{ij}$, and then run the
above search over the corresponding box.

A non-systematic search for D4 equations whose analytic solution expands
as a series in $\Z [[t]]$ was pioneered by Almkvist, van Enckevort,
van Straten and Zudilin, \cite{AZ},  \cite{ES}. See \cite{ES} for a
systematic approach to recognizing a given globally nilpotent
D4 equation as the mirror DE of a Calabi-Yau family by computing invariants of
its global monodromy.

The hypergeometric pullback conjecture suggests a (presumably) more
restrictive candidate list, but it is not clear
how one can  identify these among all D4 equations,
without further assumptions.

\pph \bf Del Pezzo surfaces and D2's. \rm

The only rank $1$ del Pezzo is
$\P ^2$, so it might seem that our program is just not applicable here.
However, it turns out (Orlov and Golyshev, unpublished) that the
three-dimensional subspace of
the total cohomology generated by the classical powers of the anticanonical
class is stable under quantum multiplication by the anticanonical class, and
it gives rise to D2 equations for del Pezzo surfaces of degrees
$9,6,5,4,3,2,1.$ In \cite{Go1}
the parametric D2 equation was identified with a particular case
of the classical Heun equation that
had been studied by Beukers \cite{Be} and Zagier \cite{Za}. Zagier had
run a search
over a large box for D2 equations with analytic solution in $\Z
[[t]]$,
see the list in \cite{Za}. Our counting D2
equations are hypergeometric in degrees $4,3,2,1$ and are
hypergeometric pullbacks in degrees $9,6,5.$

The classification of del Pezzo surfaces is of course well known;
however, it might be interesting to understand the significance of
the non-D2 equations arising as canonical pullbacks in degrees
$8$ and $7$.

\pph \bf Singular Fano 3-folds. \rm The classification of singular
Fano 3-folds of Picard rank 1 is of interest in birational geometry,
see \cite{CPR}. Corti has suggested to extend the mirror approach to
the classification of $\Q$-Fano 3-folds with prescribed (say terminal,
or canonical) singularities. One expects that to
a $\Q$-Fano 3-fold one can associate a differential equation
that reflects its properties
in much the same way as D3 equations do for smooth 3-folds. In order
to construct it as a counting DE one would have to rely on a theory of
Gromov-Witten invariants of singular varieties, which is not yet
sufficiently developed. A provisional solution is to model the
construction of such a DE on the known smooth examples, formally
generalizing them in the simplest cases such as complete intersections
in weighted projective spaces.

An instance  of a pair of mirror dual  problems in this setup is
due to Corti and myself. Let
$\P (w_0,w_1,w_2,w_3)$ be a weighted projective space,
$d= \sum w_j.$ The operator
$$\prod_{j=0}^3 (w_i^{w_i}(D-\frac{w_i-1}{w_i})(D-\frac{w_i-2}{w_i})\dots D) -
d^d  t(D+\frac{1}{d})(D+\frac{2}{d})\dots (D+\frac{d-1}{d})(D+1)$$
gives rise to a hypergeometric D-module, whose essential constituent
we call the anticanonical Riemann-Roch D-module. It is easy to show that
the monodromy of this D-module respects a real orthogonal form.
The problem of classification of weighted
$\P^3$ with canonical singularities happens to admit a mirror dual problem:
to classify the anticanonical Riemann-Roch D-modules such
that the form above is of signature $(2, n-2)$.

\bigskip
\bigskip
\bigskip

I am obliged to Helena Verrill who checked the modular formulas and
made a number of suggestions.

I thank  Alexander Givental for references to quantum Lefschetz, Victor
Przyjalkowski for explanations, Constantin Shramov
and Jan Stienstra for comments. I thank
Yuri Manin for his interest
in my work.

Special thanks go to Alessio Corti who organized the
lecture series at Cambridge and edited these notes.

\end{document}